%% file: move.tex
\input makro.tex

\gauche={S. Burckel}

\centerline{\XIV Natural Moves for Knots and Links.}

\bn {\bf Abstract. We propose some natural generalizations of Reidemeister moves that do not increase the number
of crossings in the generated diagrams. Experimentations make us conjecture that this class of monotonic moves is
complete for computing canonical forms and then deciding isotopy.}

\bn
\ttparag{1. Introduction.}

It is a classical result that Reidemeister moves enable to transform any link diagram into 
any other equivalent one. However, that is not completely satisfactory for computations since
the identification of two $n$ crossings equivalent diagrams may need to generate intermediary 
diagrams with $m>n$ crossings. Some works have been done in order to bound this $m$ 
at least for the unknotting problem (see [2]).
\bn
Here, we investigate another direction : find new types of moves that do not increase the 
number of crossings. Hence, the possible generated diagrams are in finite number (up to their
Gauss codes representations) and we obtain a possible efficient decision algorithm.
\bn
We have not yet proved that these moves are complete, that is, enable to transform every pair
of equivalent link diagrams in a common one. However, we have already 
implemented these moves for knots. That enabled to rebuild the tables of prime knots for $n=3\dots 12$ crossings.
For larger $n$, the computations have to be done.

\bn
\ttparag{2. Definitions.}

First we propose a more general notion of strands.

\Def{(fragment).}{ Let $X$ be a link projection. A {\it fragment} of $X$ is a
continuous part of the curve $X$. A fragment is said {\it close} if it has no extremities.
Otherwise, it has two extremities and is said {\it open}.
A collection $C$ of fragments of $X$ is said {\it complete} if every crossing occuring in $C$ has its two occurrences in $C$.
}

Consider the circular Gauss code $G$ of $X$. One can see a fragment of $X$ as a part of $G$.
For instance, the trefoil projection $X$ with Gauss code $G=AbCaBc$ has fragments associated
to the parts $AbCaB$ or $cA$ of $G$. The fragment $bCaBcA$ is closed.
The collection $\{bC,Bc\}$ is complete.

Now, we define a more general notion of braid.

\Def{(box).}{ 
For any $n,m\ge 0$, a {\it $n,m$-box} is a complete
collection $C$ of fragments 
embedded in a planar square where $n$ (resp. $m$) extremities are on the top (resp. bottom).
Denote $C^+$ (resp. $C^-$) the list of extremities at the top (resp. bottom) of $C$.
The {\it product} $A.B$ of a $n,m$-box $A$ with a $m,p$ box $B$ is obtained
by the one-to-one identification of the $m$ extremities at the bottom of $A$ with the $m$ ones at the top of $B$, i.e., 
by $A^-=B^+$
}

Here is a $3,3$-box and a $4,2$ box made of $4$ fragments (one is closed)~:

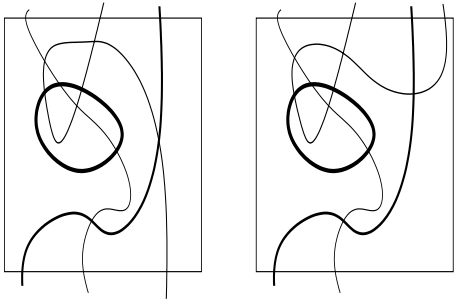
$$\dessin(box1.eps)(6)(4)$$ 

Observe that every $n$-strands braid is a $n,n$-box.

\bn
\ttparag{3. The moves.}

First, we keep the basical Reidemeister moves $R1$ and $R2$ but only in the directions that
decrease the number of crossings : $R1^-,R2^-$. We also keep $R3$ moves since they do not increase the number of crossings.

Second, we generalize $R3$ in two ways : a $2D$ operation and a $3D$ operation.

\Def{(shift).}{ Let $C$ be a collection of fragments
made of a $n,m$-box $B$ and a fragment $F$ that crosses over $B^+$. 
The {\it shift} of $C$ is obtained by moving
the fragment $F$ over $B^-$.
}
For example~:
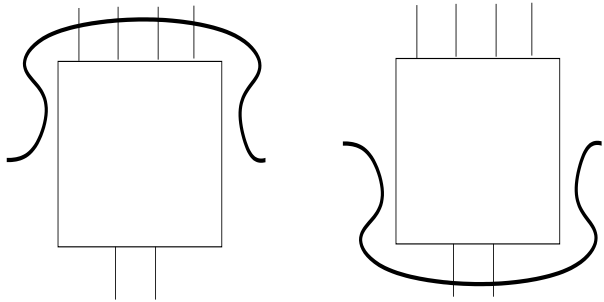
$$\dessin(shift.eps)(8)(4)$$

The shift move is a $2D$ generalization of the $R3$ move~:

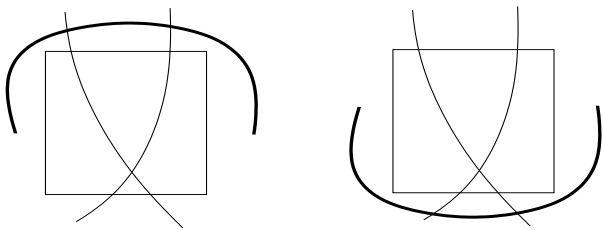
$$\dessin(r3.eps)(8)(3)$$ 

Consider $\Delta_n$, the Garside's fundamental $n$-strands braid (see [1]) and
imagine it in $3D$ as a rotation of $n$ parallel strands in the space. We do the same
for every $n$ boxes.

\Def{(rotation).}{ The {\it rotation} $\bar{C}$ of an $n,m$-box $C$ is obtained
by inverting the order of its extremities and its crossings.}

The rotation of the first example is~:

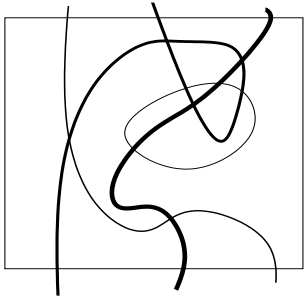
$$\dessin(rbox1.eps)(4)(4)$$

\Def{(twist).}{ 
Let $A$ be a $n,m$-box. The {\it twist} of $A$ is obtained by a complete rotation of $A$ after fixing its extremities.
We obtain the $n,m$ box $\Delta_n.\bar{A}.\Delta_m^{-1}$.
}
For example~:
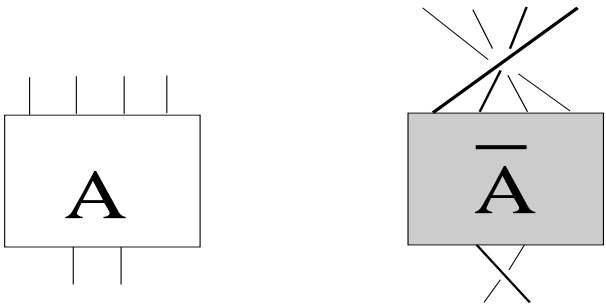
$$\dessin(twist.eps)(8)(4)$$ 

The twist move is a generalization of the $R1$ and $R2$ moves~:

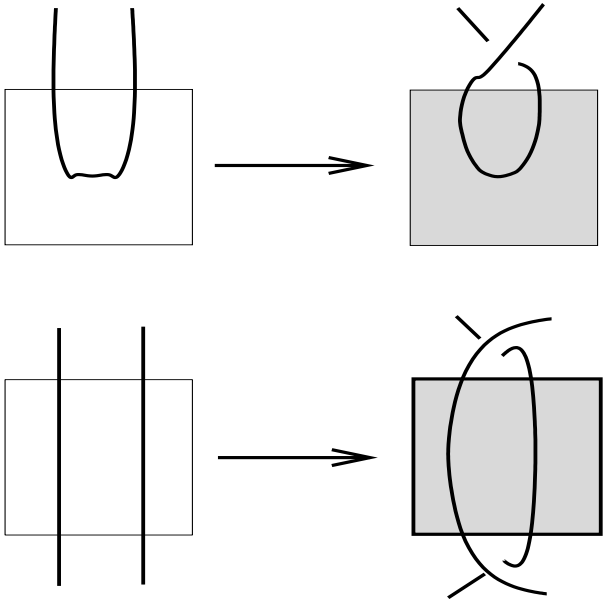
$$\dessin(genr.eps)(8)(8)$$ 

Now we make the restriction that gives monotonic moves. Of course shift and twist may introduce new crossings. However, we perform then 
reductions of diagrams by sequences of $R1^-,R2^-,R3$ moves.

\Def{(GENERIC).}{ 
A {\it generic} move consists in the application of a shift or a twist followed by a sequence of $R1^-,R2^-,R3$ moves
such that a projection $X$ with $n$ crossings is transformed in a diagram $Y$ with $m\le n$ crossings.
}
\bn
\ttparag{4. Implementation.}

It is not very difficult to implement generic moves on Gauss Codes representations.
Moreover, using the combinatorial characterizations of realizable Gauss codes (see [5]),
one can decide if some fragments of a diagram can be drawn in a $n,m$-box, considering the 
realizability of the box circled by a single curve.
\bn
We have checked generic moves for knots, considering all possible diagrams on $n=3\dots 12$ crossings
(using again [5]) and 
we have obtained reduced knots diagrams (up to mirror images). 
Computations remain to be done for $n\ge 13$.
In order to deal with amphicheirality, one has to consider the mirror images of the moves
instead of the mirror images of the diagrams.
\bn
We have obtained the expected sets of prime knots according to [6]. 
Moreover, computations are quite fast. For instance, the Perko pair identification is 
done by a single generic move from a twist of $3$ fragments. Programs and tables are available on~:
\sk
\sk
http://www2.univ-reunion.fr/\tilde burckel/box.html
\bn

\ttparag{5. Questions.}

Question 1~: Can one prove that generic moves are complete for knots (links)?
\sk
Question $\bar{1}$~: Can one find two projections $X,Y$ of equivalent knots (links) that cannot be
identified by generic moves ?
\sk 
Question 2~: Can one generalize generic moves ?
\sk
Question 3~: Bounding the number of crossings seems to be quite natural in order to produce effective and efficient algorithms. 
Are there other quantities that could be bounded in some well orderings?
\sk 
Question 4~: Some particular twist on $2$ fragments is usually called {\it flype}. 
The Tait Flyping Conjecture [7] says that flypes enable to decide the isotopy of alternating knots. 
This important statement has been proved by Kauffman [3], Murasugi [4] and Thistlethwaite [8].
We conjecture that generic moves on $3$ fragments enable to decide the isotopy of knots in general.
\bn


\def \it#1{#1}

\parindent=1.8truecm \def\Ref#1; #2; #3; #4;{ \sgni\item{\hbox to 1.6truecm{ [#1]\hfill}}{\sc #2},
{\sl #3}, #4. } \def\Reff#1; #2; #3; #4; #5; #6; #7;{ \sgni\item{\hbox to 1.6truecm{
[#1]\hfill}}{\sc #2}, {\sl #3}, #4 {\bf #5} (#6) #7. }

\bn

\Reff {1}; F. Garside; The Braid Group and other Groups; Quart. J. Math. 
Oxford; 20; 1969; 235--254;

\Reff {2}; J. Hass and J.C. Lagarias; The number of Reidemeister Moves Needed for Unknotting;
J. Amer. Math. Soc.; 14; 2001; 399--428;

\Reff {3}; L. Kauffman; Combinatorics and knot theory;
Contemp. Math.; 20; 1983; 181--200;

\Reff {4}; K. Murasugi; Jones polynomials and classical conjectures in knot theory;
Topology; 26; 1987; 187--194;

\Reff {5}; R. C. Read and P. Rosenstiehl; On the Gauss Crossing Problem;
Colloq. Math. Soc. Janos Bolyai; 18; 1976; 843--876;

\Ref {6}; D. Rolfsen; Knots and Links; Publish or Perish, Houston (1990);

\Reff {7}; P.G. Tait; On Knots I,II,III;
Scientific Papers; 1; 1898; 273--347;

\Reff {8}; M. Thistlethwaite; A spanning tree expansion for the Jones polynomial;
Topology; 26; 1987; 297--309;

\bn

\bn

\sk\ni Serge Burckel. 
\sk\ni Laboratoire ERMIT, 
\sk\ni D\'epartement de Math\'ematiques, 
\sk\ni Universit\'e de la R\'eunion, 
\sk\ni 97715 Saint-Denis, France. 
\sk\ni burckel@univ-reunion.fr
\sk\ni http://www2.univ-reunion.fr/\tilde burckel
\end

%% file: makro.tex
%

\long\def\comm(#1)#2!!!{}

{\gdef\editer
                                {\tolerance 3000
                                \gdef\\{{\tt\char'134}}\gdef\#{{\tt\char'043}}          
        
                                \gdef\{{{\tt\char'173}}\gdef\}{{\tt\char'175}}
                                \long\gdef\comm(##1)##2!!!
                                                        
        {{\leftskip20mm\parindent=0mm\medskip\hskip-20mm
                                                                \setbox1=\hbox{\\\tt##1}
                                                                {\ifdim\wd1<17mm\hbox to 
19mm{\box1\hfill}%
                                                                \else\hbox to \wd1{\unhbox1} : 
\fi}##2.\par}}
                                \centerline{\bfXII Les commandes de DehTeX}\bigskip
                                \datef\bigskip\bigskip
                                Num\'eros des familles de fontes: 
        0=cmr$\fl$rm; 1=cmmi$\fl$mi; 2=cmsy$\fl$sy; 3=cmex;  4=cmti$\fl$it;
                                5=cmsl$\fl$sl; 6=cmbx$\fl$bf; 7=cmtt$\fl$tt; 8=cmb$\fl$be; 
                                9=cmmib$\fl$bm; A=msam$\fl$sa; B=msbm$\fl$sb; 
C=eufm$\fl$go;
                                D=eurm$\fl$eu; E non utilis\'e; F non utilis\'e;
                                \def\dump{}}}


\def\<<{{\cyr\char'074~}}
                                \comm(<<)Ecrit \<<(avec blanc inseccable apr\`es)!!!

\def\>>{~{\cyr\char'076}}
                                \comm(>>)Ecrit \>> (avec blanc inseccable avant)!!!

\def\\{\hbox{$\backslash$}}
                                \comm($\\$)Ecrit $\\$; ne n\'ecessite pas de blanc apr\`es!!! 

\let\(=\langle
                                \comm(()Ecrit $\($; mode math\'ematique!!!

\let\)=\rangle
                                \comm({)})Ecrit $\)$; mode math\'ematique!!!

\def\0{{\mi 0}}\def\1{{\mi 1}}\def\2{{\mi 2}}\def\3{{\mi 3}}
\def\4{{\mi 4}}\def\5{{\mi 5}}\def\6{{\mi 6}}\def\7{{\mi 7}}
\def\8{{\mi 8}}\def\9{{\mi 9}}
                                \comm({\it chiffre})Ecrit le chiffre correspondant en 
\<<oldstyle\>>: \0,
                                                                \1, \2, \3, \4, \5, \6, \7, \8, \9!!!

\font\V=cmr5 \font\VI=cmr6 \font\VII=cmr7 \font\VIII=cmr8
\font\IX=cmr9 \font\XII=cmr10 scaled\magstep1
\font\XIV=cmr10 scaled\magstep2
\font\XVIII=cmr10 scaled\magstep3
\font\XXIV=cmr10 scaled\magstep4
                                \comm({{\it chiffre romain}})Appelle la taille de caract\`ere
                                correspondante dans la fonte cmr; disponibles: {\tt\\V, \\VI, \\VII,
                                \\VIII, \\IX, \\XII, \\XIV, \\XVIII, \\XXIV}!!!

\mathchardef\a="710B
                                \comm(a)Ecrit $\a$; mode math\'ematique!!!

\def\adots{\mathinner{\mkern2mu\raise1pt\hbox{.}
                                                                \mkern3mu\raise4pt\hbox{.}
                                                                \mkern1mu\raise7pt\hbox{.}}}
                                \comm(adots)Ecrit $\adots$; mode math\'ematique!!!

\let\al=\aleph
                                \comm(al)Ecrit $\al$; mode math\'ematique!!!

\def\and{\hbox{ and }}
                                \comm(and)Ecrit $\and$ en romain avec blanc avant et apr\`es, 
quel
                                                                que soit le mode (le m\^eme en 
pench\'e s'appelle {\tt\\sland})!!!

\def\ape{\mathchar"3A44}
                                \comm(ape)Ecrit $\ape$; mode math\'ematique!!!
                 
\def\app{\mathord{\in}}
                                \comm(app)Ecrit $\app$ sans blanc avant ni apr\`es; mode 
math\'ematique!!!

\let\aps=\triangleright
                                \comm(aps)Ecrit $\aps$; mode math\'ematique!!!
                
\let\av=\triangleleft
                                \comm(av)Ecrit $\av$; mode math\'ematique; autre nom: {\tt 
avs}!!!
                
\def\ave{\mathchar"3A45}
                                \comm(ave)Ecrit $\ave$; mode math\'ematique!!!

                                \comm(avs)Voir {\tt ave}!!!

\mathchardef\b="710C       
                                \comm(b)Ecrit $\b$; mode math\'ematique!!!

\def\bbar#1{\overline{\mathstrut#1}}
                                \comm(bbar)La commande {\tt\\bbar\#} place une barre sur \#, \`a
                                                                une hauteur uniforme et assez 
grande; mode math\'ematique!!!

                                                                \font\cmbX=cmb10 at 10pt 
\font\cmbVII=cmb10 at 7pt 
                                                                \font\cmbV=cmb10 at 5pt 
\textfont8=\cmbX
                                                        
        \scriptfont8=\cmbVII\scriptscriptfont8=\cmbV
                                                                \def\be{\fam8\cmbX}
                                                                \let\bfe=\be
                                \comm(be)Passe \`a la famille \<<gras \'etroit\>> (fontes cmb); 
{\be
                                                                Ceci est du be}; autre nom: 
{\tt\\bfe}!!!

                                                        
        \font\cmbxX=cmbx10\font\cmbxVII=cmbx7 
                                                                \font\cmbxV=cmbx5 
\textfont6=\cmbxX
                                                        
        \scriptfont6=\cmbxVII\scriptscriptfont6=\cmbxV
                                                                \def\bf{\fam6\cmbxX}

\font\bfXII=cmbx10 scaled\magstep1  
              \font\bfXIV=cmbx10 scaled\magstep2  
              \font\bfXVIII=cmbx10 scaled\magstep3  
              \font\bfXXIV=cmbx10 scaled\magstep4
                                \comm(bf{{\it chiffre romain}})Passe \`a la fonte \<<gras\>>\ 
(cmbx)
                                                                de la taille correspondante;
                                                                disponibles: {\tt\\bfXII, \\bfXIV, 
\\bfXVIII, \\bfXXIV}!!!

\def\bg{\medbreak\medskip}
        \comm(bg)Provoque un grand saut et favorise la coupure; la moiti\'e du
                                saut est supprim\'ee apr\`es un display!!!

\def\bgni{\medbreak\medskip\noindent}
                                \comm(bgni)Provoque un grand saut, favorise la coupure et
                                supprime l'indentation apr\`es!!!

\let\bi=\cmbxtiX
                                \comm(bi)Passe \`a la fonte \<<gras italique\>> (fontes cmbxti); 
{\bi
                                Ceci est du bi}!!!

                                \comm(bl)Ecrit $\bullet$!!!

                                                        
        \font\cmmibX=cmmib10\font\cmmibVII=cmmib7
                \font\cmmibV=cmmib5 \textfont9=\cmmibX
                \scriptfont9=\cmmibVII\scriptscriptfont9=\cmmibV
                \def\bm{\fam9\cmmibX}
                                \comm(bm)Appelle la famille \<<gras math\'ematique\>> (fontes
                                                                cmmib); {\bm Ceci est du bm}!!!

\def\bn{\bigskip\nobreak}
                                \comm(bn)Provoque un saut et d\'efavorise la coupure!!!

                                \comm(bnni)Provoque un saut, d\'efavorise la coupure et 
supprime
                                                                l'indentation!!!

\long\def\boxit#1#2{\setbox1=\hbox{\kern#1{#2}\kern#1}%
                                                                \dimen1=\ht1 \advance\dimen1 by 
#1
                                                                \dimen2=\dp1 \advance\dimen2 by 
#1
                                                                \setbox1=\hbox{\vrule 
height\dimen1 depth\dimen2\box1\vrule}%
                                                        
        \setbox1=\vbox{\hrule\box1\hrule}%
                                                                \advance\dimen1 by .4truept 
\ht1=\dimen1
                                                                \advance\dimen2 by .4truept 
\dp1=\dimen2 \box1\relax}
                                \comm(boxit)Encadre son second argument; syntaxe
                                                                {\tt\\boxit$\{$2pt$\}\{$xx$\}$} 
donne \boxit{2pt}{{\tt xx}}!!!

\let\bs=\cmbxslX
                                \comm(bs)Passe \`a la fonte \<<gras pench\'e\>> (fontes cmbxsl); 
{\bs
                                ceci est du bs}!!!

\def\card{{\rm card}}
                                \comm(card)Ecrit \card!!!               

\def\carre{\hbox{\sa\char'004}}
                                \comm(carre)Ecrit \carre!!!

\catcode`\@=11\def\cases#1{\left\{\,\vcenter{\normalbaselines\m@th
                                                        
        \ialign{$##\hfil$&\quad{\rm##}\hfil\crcr#1\crcr}}\right.}
                                                                \catcode`\@=12

\def\cf{\hbox{\it c.f. }}
                                \comm(cf)Ecrit \cf\ avec un espace apr\`es!!!

\font\cmdunhX=cmdunh10
                                \comm(cmdunhX)Passe \`a la fonte cmdunh10; {\cmdunhX Ceci 
est du
                                                                cmdunh10}!!!

\def\comp{{\scriptscriptstyle\circ}}
                                \comm(comp)Ecrit $\comp$; mode math\'ematique!!!

                                \comm(Cor)La commande {\tt Cor\#} provoque un saut, puis 
\'ecrit
                                                                sans indentation {\bf Corollary\#.-
}!!!

                                \comm(cyr)Fait passer dans la fonte \<<wncyr10\>>; {\cyr Ceci 
est du
                                                                cyr}!!!

\mathchardef\d="710E
                                \comm(d)Ecrit $\d$; mode math\'ematique!!!

\mathchardef\D="7101
                                \comm(D)Ecrit $\D$; mode math\'ematique!!!

\def\date{\par\line{\hfil\ifcase\month\or January\or
                                                                February\or March\or April\or 
May\or June\or July\or August\or
                                                                September\or October\or 
November\or December\fi\ 
                                                                {\oldstyle\the\day},\ 
{\oldstyle\the\year}}}
                                \comm(date)Met la date en anglais \`a droite de la ligne!!!

\def\datef{\par\line{\hfil le\ {\oldstyle\the\day}\
                 \ifcase\month\or
                 janvier\or f\'evrier\or mars\or avril\or mai\or juin\or
                 juillet\or ao\^ut\or septembre\or octobre\or novembre\or
                 d\'ecembre\fi\ {\oldstyle\the\year}}}
                                \comm(datef)Met la date en francais \`a droite de la ligne!!!

\def\Def{\bgni{\bf Definition.- }}
                                \comm(Def)Provoque un saut, puis \'ecrit {\bf Definition.-}!!!

\let\der=\partial
                                \comm(der)Ecrit $\der$; mode math\'ematique!!!

                                \comm(diagram)Pour coder un diagramme; m\^eme syntaxe que
                                                                {\tt\\matrix}, simplement les espaces 
sont red\'efinis!!!

\def\Dom{\mathord{\rm Dom}}
                                \comm(Dom)Ecrit $\Dom$; mode math\'ematique!!!

\mathchardef\e="7122
                                \comm(e)Ecrit $\e$; mode math\'ematique!!!

\def\eg{\hbox{\it e.g. }}
                                \comm(eg)Ecrit \eg avec un espace apr\`es!!!

                                \comm(entete)Place l'entete de Caen en haut de page; attention! il 
faut
                                                                que le sceau soit dans le fichier 
pictures!!!

\let\equ=\Longleftrightarrow
                                \comm(equ)Ecrit $\equ$; mode math\'ematique!!!

                                \comm(esp)Force un blanc, m\^eme dans les fontes o\`u l'espace 
est nul,
                                                                comme cmmi qui sert dans 
oldstyle!!!

\def\et{\hbox{ et }}
                                \comm(et)Ecrit $\et$ en tout mode (avec un espace avant et 
apr\`es)!!!

\def\etc{\hbox{\it etc\dots\ }}
                                \comm(etc)Ecrit \etc avec un espace apr\`es!!!

                                                        
        \font\euX=eurm10\font\euVII=eurm7\font\euV=eurm5
                                                                \textfont13=\euX
                                                        
        \scriptfont13=\euVII\scriptscriptfont13=\euV
                                                                \def\eu{\fam13\euX}
                                \comm(eu)Passe \`a la famille \<<euler\>> (fontes eurm); {\eu
                                                                Ceci est du eu}!!!

\mathchardef\eup="0D3A
                                \comm(eup)Ecrit $\eup$ (point en fonte \<<euler\>>); mode
                                                                math\'ematique!!!

\mathchardef\euv="0D3B
                                \comm(euv)Ecrit $\euv$ (virgule en fonte \<<euler\>>); mode
                                                                math\'ematique!!!

\let\ev=\emptyset
                                \comm(ev)Ecrit $\ev$; mode math\'ematique!!!

\let\ex=\exists
                                \comm(ex)Ecrit $\ex$; mode math\'ematique!!!

                                \comm(Ex)Provoque un grand saut, puis \'ecrit {\bf Example.-
}!!!

\mathchardef\f="7127
                                \comm(f)Ecrit $\f$; mode math\'ematique!!!

\mathchardef\F="7108
                                \comm(F)Ecrit $\F$; mode math\'ematique!!!
                 
\let\fl=\longrightarrow
                                \comm(fl)Ecrit $\fl$; mode math\'ematique!!!

                                \comm(Fin)Ecrit $\carre$ pr\'ec\'ed\'e d'un blanc inseccable!!!

\let\flc=\rightarrow
                                \comm(flc)Ecrit $\flc$ (fl\`eche courte); mode math\'ematique!!!

{\catcode`@=11
                                                                \catcode`\;=\active%
                                                                \catcode`\:=\active%
                                                                \catcode`\!=\active%
                                                                \catcode`\?=\active%
                                \gdef\francais{
                                                                \catcode`\;=\active%
                                                        
        \def;{\relax\ifhmode\ifdim\lastskip>\z@\unskip\fi
                                                                \kern.2em\fi\string;}
                                                                \catcode`\:=\active%
                                                        
        \def:{\relax\ifhmode\ifdim\lastskip>\z@\unskip\fi
                                                                \kern.2em\fi\string:}
                                                                \catcode`\!=\active%
                                                        
        \def!{\relax\ifhmode\ifdim\lastskip>\z@\unskip\fi
                                                                \kern.2em\fi\string!}
                                                                \catcode`\?=\active%
                                                        
        \def?{\relax\ifhmode\ifdim\lastskip>\z@\unskip\fi
                                                                \kern.2em\fi\string?}
                                                                \frenchspacing\catcode`\@=12}}
                                \comm(francais)Adapte \`a la typographie francaise;
                                                                redefinit les caracteres ; : ? !  aux 
normes francaises, c'est \`a
                                                                dire avec un blanc inseccable devant; 
cela ne change rien de taper
                                                                ou de ne pas taper les blancs avant 
les caracteres ; : ! ? (par contre
                                                                cela change d'en taper apres ou 
non)!!!

\mathchardef\g="710D
                                \comm(g)Ecrit $\g$; mode math\'ematique!!!

\mathchardef\G="7100
                                \comm(G)Ecrit $\G$; mode math\'ematique!!!
                  
                                                        
        \font\eufmX=eufm10\font\eufmVII=eufm7\font\eufmV=eufm5
                                                        
        \textfont12=\eufmX\scriptfont12=\eufmVII
                                                                \scriptscriptfont12=\eufmV
                                                                \def\go{\fam12\eufmX}
                                \comm(go)Passe dans la famille \<<gothique\>> (fontes eufm); 
                                                                {\go Ceci est du go}!!!

\mathchardef\gpref="3A40
                                \comm(gpref)Ecrit $\gpref$ (le g signifie \<<grand\>>) ; mode
                                                                math\'ematique!!!
                                                        
\let\gprefe=\sqsubseteq
                                \comm(gprefe)Ecrit $\gprefe$ (le g signifie \<<grand\>>) ; mode 
                                                                math\'ematique!!!

\mathchardef\gsuff="3A41
                                \comm(gsuff)Ecrit $\gsuff$ (le g signifie \<<grand\>>) ; mode
                                                                math\'ematique!!!
                                                        
\let\gsuffe=\sqsupseteq
                                \comm(gsuffe)Ecrit $\gsuffe$ (le g signifie \<<grand\>>) ; mode 
                                                                math\'ematique!!!

\mathchardef\h="7112
                                \comm(h)Ecrit $\h$; mode math\'ematique!!!

\mathchardef\H="7102
                                \comm(H)Ecrit $\H$; mode math\'ematique!!!
                  
\font\hel=cmb10 at 10 pt
               \font\helVII=cmb10 at 7 pt
               \font\helXII=cmb10 at 12 pt
               \font\helXIV=cmb10 at 14 pt
                                \comm(hel)Passe dans la fonte \<<cmb10\>>; 
                                                                {\hel Ceci est du hel};
                                                                suivie \'eventuellemnt d'une taille en 
romain majuscule;
                                                                disponibles: {\tt \\helVII, \\helXII, 
\\helXIV}!!!

\font\helbf=cmb10 at 10 pt
               \font\helbfXII=cmb10 at 12 pt
               \font\helbfXIV=cmb10 at 14 pt
               \font\helbfXVIII=cmb10 at 18 pt
               \font\helbfXXIV=cmb10 at 24 pt
               \font\helbfXXXVI=cmb10 at 36 pt
                                \comm(helbf)Passe \`a la fonte \<<cmb10 gras\>>;
                                                         {\helbf ceci est du helbf};
                                                                suivie \'eventuellemnt d'une taille en 
chiffre romain majuscule; 
                                                                disponibles:
                                                                {\tt\\helbfXII, \\helbfXIV, 
\\helbfXVIII,\\helbfXXIV,
                                                                \\helbfXXXVI}!!!

\font\helbi=cmb10 at 10 pt
               \font\helbiXII=cmb10 at 12 pt
               \font\helbiXIV=cmb10 at 14 pt
               \font\helbiXVIII=cmb10 at 18 pt
               \font\helbiXXIV=cmb10 at 24 pt
               \font\helbiXXXVI=cmb10 at 36 pt
                                \comm(helbi)Passe \`a la fonte \<<cmb10 gras italique\>>; 
                                                                {\helbi Ceci est du helbi}; 
                                                                suivie \'eventuellemnt d'une taille en 
chiffres romains majuscules;
                                                                disponibles:
                                                                {\tt\\helbiXII, \\helbiXIV, 
\\helbiXVIII, \\helbiXXIV, 
                                                                \\helbiXXXVI}!!!

\def\hfl#1#2{\smash{\mathop{\hbox to 12truemm{\rightarrowfill}}
                                                        
        \limits^{\scriptstyle#1}_{\scriptstyle#2}}}
                                \comm(hfl{\#1}{\#2})Trace une fl\`eche horizontale de 12 mm 
avec
                                                                {\tt\#1} dessus et {\tt\#2} dessous!!!

                                \comm(hhat)Ecrit un grand chapeau!!!

\def\ie{\hbox{\it i.e. }}
                                \comm(ie)Ecrit \ie avec un blanc apr\`es!!!

\def\iff{if and only if }
                                \comm(iff)Ecrit \iff avec un blanc apr\`es!!!

\def\Im{{\be Im}}
                                \comm(Im)Ecrit \Im!!!

\let\imp=\Longrightarrow
                                \comm(imp)Ecrit $\imp$; mode math\'ematique; autre nom:
                                                                {\tt\\impl}!!!

\let\impl=\Longrightarrow
                                \comm(impl)Voir {\tt\\imp}!!!

\let\inc=\subset
                                \comm(inc)Ecrit $\inc$; mode math\'ematique!!!

\let\ince=\subseteq
                                \comm(ince)Ecrit $\ince$; mode math\'ematique!!!

\def\inserer(#1)(#2)(#3)(#4)(#5){
                                                                \dimen1=#2\divide\dimen1 
by1000\multiply\dimen1 by#4
                                                                \dimen2=#3\divide\dimen2 
by1000\multiply\dimen2 by#4
                                                        
        \midinsert\vglue\dimen2$$\vbox{\hsize=\dimen1
                                                                \ni\special{picture #1 scaled #4}}$$
                                                                \centerline{#5}\endinsert}
                                \comm(inserer)Ins\`ere au mieux un dessin
                                                                la syntaxe est {\tt\\inserer(nom de la 
figure)(largeur)(hauteur)
                                                                (echelle)(commentaire)}
                                                                (les parametres largeur et hauteur
                                                                sont \`a lire dans le fichier Pictures; 
il faut taper l'unite)!!!

\let\inter=\cap
                                \comm(inter)Ecrit $\inter$; mode math\'ematique!!!

\let\Inter=\bigcap
                                \comm(Inter)Ecrit $\Inter$; mode math\'ematique!!!

                                                        
        \font\cmtiX=cmti10\font\cmtiVII=cmti7 
                                                                \font\cmtiV=cmti10 at 5pt 
                                                        
        \textfont4=\cmtiX\scriptfont4=\cmtiVII\scriptscriptfont4=\cmtiV
                                                                \def\it{\fam4\cmtiX}

\mathchardef\k="7114
                                \comm(k)Ecrit $\k$; mode math\'ematique!!!

\mathchardef\l="7115
                                \comm(l)Ecrit $\l$; mode math\'ematique!!!

\mathchardef\L="7103
                                \comm(L)Ecrit $\L$; mode math\'ematique!!!
                  
\def\Lem#1{\bgni{\bf Lemma#1.-}}
                                \comm(Lem)La commande {\tt\\Lem\#} provoque un saut, puis 
\'ecrit
                                                                sans indentation {\bf Lemma\#.-}!!!

                                \comm(lettre)Ecrit le haut de page pour une lettre en anglais; il faut
                                                                le sceau dans le fichier pictures!!!

                                \comm(lettref)Ecrit le haut de page pour une lettre en francais; il 
faut
                                                                le sceau dans le fichier pictures!!!

                                \comm(lettrem)Ecrit le haut de page pour une lettre en francais!!!

\def\lpol{\leavevmode\setbox0=\hbox{l}\hbox to \wd0{\hss\char'40l}}
                                \comm(lpol)Ecrit \lpol\ (l polonais)!!!

\def\Lpol{\leavevmode\setbox0=\hbox{L}\hbox to \wd0{\hss\char'40L}}
                                \comm(Lpol)Ecrit \Lpol\ (L polonais)!!!

\mathchardef\m="7116
                                \comm(m)Ecrit $\m$; mode math\'ematique!!!

                                \comm(mg)Provoque un saut moyen et favorise la coupure!!!

                                \comm(mgni)Provoque un saut moyen, favorise la coupure et 
supprime
                                                                l'indentation!!!

                                                                \def\mi{\fam1\teni}
                                \comm(mi)Passe \`a la famille \<<math. italique\>> (fontes cmmi); 
{\mi
                                                                ceci est du mi}; autres noms: 
{\tt\\mit} (mode math\'ematique),
                                                                {\tt\\oldstyle} (tout mode)!!!

                                \comm(mn)Provoque un saut moyen et d\'efavorise la coupure!!!

                                \comm(mnni)Provoque un saut moyen, d\'efavorise la coupure et
                                                                supprime l'indentation!!!

\def\move#1#2#3
                                                                {\vbox to0pt{\kern-
#3mm\smash{\hbox{\kern#2mm{#1}}}\vss}
                                                                \nointerlineskip}
                                \comm(move\#1\#2\#3)Ecrit \#1 \`a une abscisse de \#2 mm et \`a
                                                                une ordonn\'ee de \#3 mm par 
rapport au point courant (qui n'est
                                                                pas modifi\'e); mode vertical. 
Utilisable pour figures. Fabriquer
                                                                une \\{\tt vbox} suivant la syntaxe 
suivante \par
                                                                \hskip3cm{\tt \\vbox to ... 
mm\{\\hsize=... mm\\vfill\par
                                                         \hskip3cm\\special \{picture ... scaled ... \} 
\par 
                                                         \hskip3cm\\move \{\ ... \}\{\ ... \}\{\ ... \} 
\par 
                                                                \hskip3cm\\hfill\} }\par
                                                                Le {\tt\\vfill} est pour tasser la boite 
vers le bas, le {\tt\\hfill}
                                                                est pour tasser vers la gauche et 
assurer la largeur!!!

\mathchardef\n="7117
                                \comm(n)Ecrit $\n$; mode math\'ematique!!!

\def\Nat{\hbox{\rm I\kern-.9truemm\hbox{N}}}
                                \comm(Nat)Ecrit \Nat!!!

\let\ni=\noindent
                                \comm(ni)Supprime l'indentation!!!

\mathchardef\NN="0B4E
                                \comm(NN)Ecrit $\NN$: mode math\'ematique!!!

\def\non{\mathord{\be non}}
                                \comm(non)Ecrit $\non$; mode math\'ematique!!!

                                \comm(notapp)Ecrit $\notin$ sans blanc avant ni apr\`es; mode
                                                                        math\'ematique!!!

\mathchardef\o="7121
                                \comm(o)Ecrit $\o$; mode math\'ematique!!!

\mathchardef\O="710A
                                \comm(O)Ecrit $\O$; mode math\'ematique!!!
                  
\let\oo=\infty
                                \comm(oo)Ecrit $\oo$; mode math\'ematique!!!

\def\ou{\hbox{ ou }}
                                \comm(ou)Ecrit $\ou$ en tout mode (avec un blanc avant et 
apr\`es)!!!

\mathchardef\p="7119
                                \comm(p)Ecrit $\p$; mode math\'ematique!!!

\mathchardef\P="7105
                                \comm(P)Ecrit $\P$; mode math\'ematique!!!

\def\pmb#1{\setbox0=\hbox{#1}\kern-.15pt\copy0\kern-\wd0
                 \kern-.3pt\copy0\kern-\wd0\kern-.15pt\raise.1pt\box0}
                                \comm(pmb)Met en gras son param\`etre: {\tt\\pmb\#} \'ecrit \# en
                                                                faux gras (ce qui donne 
\pmb{\#})!!!

\def\pp{\hbox{$\ldots$}}
                                \comm(pp)Ecrit \pp!!!

\def\Ppar{\mathhexbox27B}
                                \comm(Ppar)Ecrit \Ppar!!!

\def\ppp{\hbox{$, \ldots, $}}
                                \comm(ppp)Ecrit \ppp!!!

\let\pr=\vdash
                                \comm(pr)Ecrit $\pr$; mode math\'ematique!!!

                                \comm(Pr)Fait un saut et \'ecrit \<<{\it Proof.}\>>!!!

\def\pref{\mathrel{\scriptstyle\gpref}}
                                \comm(pref)Ecrit la relation $\pref$; mode math\'ematique!!!
                                                        
\def\prefe{\mathrel{\scriptstyle\gprefe}}
                                \comm(prefe)Ecrit la relation $\prefe$; mode math\'ematique!!!
                                                        
\def\Prop#1{\bgni{\bf Proposition#1.-}}
                                \comm(Prop)La commande {\tt\\Prop\#} provoque un saut, puis 
\'ecrit
                                                                sans indentation {\bf Proposition\#.-
}!!!

\mathchardef\Pw="0C50
                                \comm(Pw)Ecrit $\Pw$; mode math\'ematique!!!

\mathchardef\q="711F
                                \comm(q)Ecrit $\q$; mode math\'ematique!!!

\let\qq=\forall
                                \comm(qq)Ecrit $\qq$; mode math\'ematique!!!

\mathchardef\r="711A
                                \comm(r)Ecrit $\r$; mode math\'ematique!!!

\mathchardef\res="0A16
                                \comm(res)Ecrit $\res$; mode math\'ematique!!!
                                                        
\def\resp{\hbox{\it resp. }}
                                \comm(resp)Ecrit \resp (avec un blanc apr\`es)!!!

\def\Ree{\hbox{\rm I\kern-.9truemm\hbox{R}}}
                                \comm(Ree)Ecrit \Ree!!!

\mathchardef\RR="0B52
                                \comm(RR)Ecrit $\RR$; mode math\'ematique!!!

\def\rres{\hbox to 4.5pt{$\res$\kern-1pt\hbox{$\res$}\hss}}
                                \comm(rres)Ecrit $\rres$; mode math\'ematique!!!

\mathchardef\s="711B
                                \comm(s)Ecrit $\s$; mode math\'ematique!!!

\mathchardef\S="7106
                                \comm(S)Ecrit $\S$; mode math\'ematique!!!

                                                        
        \font\msamX=msam10\font\msamVII=msam7 
                                                                \font\msamV=msam5 
\textfont10=\msamX
                                                        
        \scriptfont10=\msamVII\scriptscriptfont10=\msamV
                                                                \def\sa{\fam10\msamX}
                                \comm(sa)Passe \`a la famille \<<symboles a\>> (fontes msam); 
{\sa
                                                                Ceci est du sa}!!!

                                                        
        \font\msbmX=msbm10\font\msbmVII=msbm7 
                                                                \font\msbmV=msbm5 
\textfont11=\msbmX
                                                        
        \scriptfont11=\msbmVII\scriptscriptfont11=\msbmV
                                                                \def\sb{\fam11\msbmX}
                                                                \let\db=\sb
                                \comm(sb)Passe \`a la famille \<<symboles b\>> (fontes msbm); 
{\sb
                                                                CECI EST DU SB}; autre nom: 
{\tt\\db}!!!

\let\sat=\models
                                \comm(sat)Ecrit $\sat$; mode math\'ematique!!!

                                \comm(sg)Provoque un saut petit et favorise la coupure!!!

\def\sgni{\smallbreak\noindent}
                                \comm(sgni)Provoque un saut petit, favorise la coupure et 
supprime
                                                                l'indentation!!!

                                                        
        \font\cmslX=cmsl10\font\cmslVII=cmsl10 at 7 pt 
                                                                \font\cmslV=cmsl10 at 5pt 
                                                        
        \textfont5=\cmslX\scriptfont5=\cmslVII\scriptscriptfont5=\cmslV
                                                                \def\sl{\fam5\cmslX}

\def\sland{\hbox{ \sl and }}
                                \comm(sland)Ecrit $\sland$ en tout mode!!!

\font\smcap=cmcsc10                                                             
                                \comm(smcap)Passe \`a la fonte \<<petites capitales\>>\ (cmcsc);
                                                                {\smcap Ceci est du smcap}!!!

                                \comm(sn)Provoque un saut petit et d\'efavorise la coupure!!!

                                \comm(snni)Provoque un saut petit, d\'efavorise la coupure et
                                                                supprime l'indentation!!!

\def\Spar{\mathhexbox278}
                                \comm(Spar)Ecrit \Spar!!!

\mathchardef\SS="0B53
                                \comm(SS)Ecrit $\SS$: mode math\'ematique!!!

\def\ssi{si et seulement si }
                                \comm(ssi)Ecrit \ssi (avec un blanc apr\`es)!!!

\def\suff{\mathrel{\scriptstyle\gsuff}}
                                \comm(suff)Ecrit la relation $\suff$; mode math\'ematique!!!
                                                        
\def\suffe{\mathrel{\scriptstyle\gsuffe}}
                                \comm(suffe)Ecrit la relation $\suffe$; mode math\'ematique!!!
                                                        
\def\Supp{\mathord{\be Supp}}
                                \comm(Supp)Ecrit $\Supp$; mode math\'ematique!!!

                                                                \def\sy{\fam2\tensy}
                                                                \let\ca=\sy
                                \comm(sy)Passe \`a la famille \<<symbole\>> (fontes cmsy); {\sy
                                                                CECI EST DU SY}; autres noms: 
{\tt\\cal} (mode math\'ematique),
                                                                {\tt\\ca} (tout mode)!!!

\mathchardef\t="711C
                                \comm(t)Ecrit $\t$; mode math\'ematique!!!

\def\Thm#1{\bgni{\bf Theorem#1.-}}
                                \comm(Thm)La commande {\tt\\Thm\#} provoque un saut, puis 
\'ecrit
                                                                sans indentation {\bf Theorem\#.-
}!!!

\font\tim=cmb10 at 10pt 
                
               \font\timXIV=cmb10 at 14pt 
               \font\timXVIII=cmb10 at 18pt 
               \font\timVII=cmb10 at 7pt 
                                \comm(tim)Passe dans la fonte \og cmb10\fg; 
                                                                {\tim Ceci est du tim};
                                                                suivie \'eventuellemnt d'une taille en 
chiffres romains majuscules;
                                                                disponibles: {\tt\\timVII, \\timlXII, 
\\timXIV, \\timXVIII}!!!

\font\timbf=cmb10 at 10pt 
               
              \font\timbfXIV=cmb10 at 14pt 
              \font\timbfXVIII=cmb10 at 18pt 
              \font\timbfVII=cmb10 at 7pt
                                \comm(timbf)Passe dans la fonte \og cmb10 gras\fg; 
                                                                {\timbf Ceci est du timbf};
                                                                suivie \'eventuellemnt d'une taille en 
chiffres romains majuscules;
                                                                disponibles: {\tt\\timbfVII, 
\\timbflXII, \\timbfXIV,
                                                                \\timbfXVIII}!!!

                                                        
        \font\cmttX=cmtt10\font\cmttVII=cmtt10 at 7 pt 
                                                                \font\cmttV=cmtt10 at 5pt 
                                                        
        \textfont7=\cmtiX\scriptfont7=\cmttVII\scriptscriptfont7=\cmttV
                                                                \def\tt{\fam7\cmttX}

                                \comm(ttil)Ecrit un grand tilde!!!

\mathchardef\u="711D
                                \comm(u)Ecrit $\u$; mode math\'ematique!!!

\mathchardef\U="7107
                                \comm(U)Ecrit $\U$; mode math\'ematique!!!

\let\un=\cup
                                \comm(un)Ecrit $\un$; mode math\'ematique!!!

\let\Un=\bigcup
                                \comm(Un)Ecrit $\Un$; mode math\'ematique!!!

\def\val{\mathop{\be val}}
                                \comm(val)Ecrit $\val$; mode math\'ematique!!!

\def\var{{\be var}}
                                \comm(var)Ecrit $\var$!!!

                                \comm(vfl{\#1}{\#2})Trace une fl\`eche verticale de 12 mm avec
                                                                {\tt\#1} \`a gauche et {\tt\#2} \`a 
droite!!!

\mathchardef\w="7120
                                \comm(w)Ecrit $\w$; mode math\'ematique!!!

\mathchardef\W="7109
                                \comm(W)Ecrit $\W$; mode math\'ematique!!!

                                \comm(wrt)Ecrit \<<with respect to\>> suivi d'un blanc 
inseccable!!!

\mathchardef\x="7118
                                \comm(x)Ecrit $\x$; mode math\'ematique!!!

\mathchardef\X="7104
                                \comm(X)Ecrit $\X$; mode math\'ematique!!!

\mathchardef\y="7111
                                \comm(y)Ecrit $\y$; mode math\'ematique!!!

\mathchardef\z="7110
                                \comm(z)Ecrit $\z$; mode math\'ematique!!!

\def\ZZ{{\db Z}}
                                \comm(ZZ)Ecrit $\ZZ$!!!

                                \comm( )\vfill\eject\centerline{\bfXII Anciennes commandes}!!!

\let\bfe=\be
                                \comm(bfe)Ancien nom de {\tt\\be}!!!

                                \comm(bfi)Ancien nom de {\tt\\bi}!!!
  
\let\ca=\sy
                                \comm(ca)Ancien nom de la famille {\tt\\sy}!!!

\font\cmmibX=cmmib10
                                \comm(cmmibX)Passe \`a la fonte cmmib10; {\cmmibX Ceci est 
du
                                                                cmmib}!!!

\let\db=\sb
                                \comm(db)Ancien nom de la famille {\tt\\sb}!!!

\def\fg{~{\cyr\char'076}}
                                \comm(fg)Ancien nom de {\tt\\>>}!!!

                                \comm(msymX)Passe \`a la fonte \<<msbm10\>>; obsol\`ete: 
nouveau
                                                                nom de famille: {\tt\\db}!!!

\def\og{{\cyr\char'074}~}
                                \comm(og)Ancien nom de {\tt\\<<}!!!

                                \comm(rd)Ancien nom de la famille {\tt\\mi}!!!

\def\picture #1 by #2 (#3){
 \vbox to #2{\hrule width #1 height 0pt depth 0pt
 \vfill\special{picture #3}}}
\def\dessin(#1)(#2)(#3)(#4){{
 \dimen0=#2 \dimen1=#3
 \divide\dimen0 by 1000 \multiply\dimen0 by #4
 \divide\dimen1 by 1000 \multiply\dimen1 by #4
 \picture \dimen0 by \dimen1 (#1 scaled #4)}}

\newcount\numlem
\newcount\numchap
\hsize=12cm
\vsize=20truecm
\voffset=1cm
\hoffset=5mm
\newtoks\gauche \gauche={Serge Burckel}
\newtoks\droite 
\def\makeheadline{\vbox to 0pt{\vskip -40pt\line{\vbox to 8.5pt{}\the
\headline}\vss\nointerlineskip}}
\headline={{\voffset 2cm\sevenbf\the\gauche\hfill\the\droite}}

\def\sk{\smallskip}

\overfullrule=0mm 
\parindent=0mm

 

\def\ca{{\bi a}}

\def\l{n}

\def\psi{{\bf cl}}
\def\x{X}
\def\fl{{\longrightarrow}}

\def\bar#1{{\overline{#1}}}

\def\H#1{{\hskip #1truemm}}
\def\BX#1{\boxit{8pt}{\vbox{#1}}}

\def\s#1 #2{{#1}_{(#2)}}
\def\O{{\sy T}}
\def\C#1 #2{{\(#1\)_{#2}}}

\def\D{\L}

\def\ttchap#1{
\numlem=0
\vfill\eject\gauche={\the\numchap.
#1}\droite={\hfill}{\bfXVIII \BX{CHAPITRE 
\the\numchap. \bn \hbox{#1}} }\bn\vskip 2cm}

\def\ttparag
#1{{\bfXII #1}\bn}

\def\Def#1#2{\bg{{\bf Definition. #1}#2}\rm\bg}

\def\Lem#1#2{\advance\numlem by
1\bg{{\bf LEMME \the\numlem. #1}\sl#2}\rm\bg}
\def\Lemprop#1#2{\advance\numlem by 1\bg{{\bf
PROPOSITION \the\numlem. #1}\sl#2}\rm\bg}
\def\Lempropr#1#2{\advance\numlem by 1\bg{{\bf
PROPRIETE \the\numlem. #1}\sl#2}\rm\bg}
\def\Thm#1#2{\advance\numlem by 1\bg{{\bf
THEOREME \the\numlem. #1}\sl#2}\rm\bg}

\font\sc=cmcsc10

\def\ttchap#1{
\numlem=0
\vfill\eject\gauche={\the\numchap.
#1}\droite={\hfill}{\bfXVIII \BX{CHAPTER  \the\numchap.\bn
\hbox{#1}} }\bn\vskip 2cm}

\def\Def#1#2{\bg{{\bf Definition. #1}#2}\rm\bg}

\def\Lem#1#2{\advance\numlem by
1\bg{{\bf Lemma \the\numlem. #1}\sl#2}\rm\bg}
\def\Lemprop#1#2{\advance\numlem by 1\bg{{\bf
Proposition \the\numlem. #1}\sl#2}\rm\bg}
\def\Prop#1#2{\bg{{\bf
Proposition  #1}\sl#2}\rm\bg}
\def\Lempropr#1#2{\advance\numlem by 1\bg{{\bf
Property \the\numlem. #1}\sl#2}\rm\bg}
\def\Thm#1#2{\advance\numlem by 1\bg{{\bf
Theorem \the\numlem. #1}\sl#2}\rm\bg}

\catcode`\Ž=\active \def Ž{\'e}
\catcode`\ˆ=\active \def ˆ{\`a}
\catcode`\=\active \def {\`u}
\catcode`\=\active \def {\`e}
\catcode`\=\active \def {\^e}
\catcode`\™=\active \def ™{\^o}
\catcode`\"=\active \def "{\^i}
\catcode`\=\active \def {\c c}

\def\fileversion{v1.0}
\def\filedate{3 Dec 91}
\immediate\write16{Document style option `epsfig', \fileversion\space
<\filedate> (edited by SPQR)}
\chardef\atcode=\catcode`\@
\catcode`\@=11\relax
\ifx\typeout\undefined%
        \newwrite\@unused
        \def\typeout#1{{\let\protect\string\immediate\write\@unused{#1}}}
\fi
\ifx\undefined\epsfig%
\else
        \typeout{EPSFIG --- already loaded} 
\fi
%
%
\ifx\undefined\fbox\def\fbox#1{#1}\fi
%
\ifx\undefined\epsfbox\input epsf\fi
%
\ifx\undefined\@latexerr
        \newlinechar`\^^J
        \def\@spaces{\space\space\space\space}
        \def\@latexerr#1#2{%
        \edef\@tempc{#2}\expandafter\errhelp\expandafter{\@tempc}%
        \typeout{Error. \space see a manual for explanation.^^J
         \space\@spaces\@spaces\@spaces Type \space H <return> \space for
         immediate help.}\errmessage{#1}}
\fi
\def\@whattodo{You tried to include a PostScript figure which 
cannot be found^^JIf you press return to carry on anyway,^^J
The failed name will be printed in place of the figure.^^J
or type X to quit}
\def\@whattodobb{You tried to include a PostScript figure which 
has no^^Jbounding box, and you supplied none.^^J
If you press return to carry on anyway,^^J
The failed name will be printed in place of the figure.^^J
or type X to quit}
%
\newwrite\@unused
%
\def\@nnil{\@nil}
\def\@empty{}
\def\@psdonoop#1\@@#2#3{}
\def\@psdo#1:=#2\do#3{\edef\@psdotmp{#2}\ifx\@psdotmp\@empty \else
    \expandafter\@psdoloop#2,\@nil,\@nil\@@#1{#3}\fi}
\def\@psdoloop#1,#2,#3\@@#4#5{\def#4{#1}\ifx #4\@nnil \else
       #5\def#4{#2}\ifx #4\@nnil \else#5\@ipsdoloop #3\@@#4{#5}\fi\fi}
\def\@ipsdoloop#1,#2\@@#3#4{\def#3{#1}\ifx #3\@nnil 
       \let\@nextwhile=\@psdonoop \else
      #4\relax\let\@nextwhile=\@ipsdoloop\fi\@nextwhile#2\@@#3{#4}}
\def\@tpsdo#1:=#2\do#3{\xdef\@psdotmp{#2}\ifx\@psdotmp\@empty \else
    \@tpsdoloop#2\@nil\@nil\@@#1{#3}\fi}
\def\@tpsdoloop#1#2\@@#3#4{\def#3{#1}\ifx #3\@nnil 
       \let\@nextwhile=\@psdonoop \else
      #4\relax\let\@nextwhile=\@tpsdoloop\fi\@nextwhile#2\@@#3{#4}}
%
%
%
\long\def\epsfaux#1#2:#3\\{\ifx#1\epsfpercent
   \def\testit{#2}\ifx\testit\epsfbblit
        \@atendfalse
        \epsf@atend #3 . \\%
        \if@atend       
           \if@verbose
                \typeout{epsfig: found `(atend)'; continuing search}
           \fi
        \else
                \epsfgrab #3 . . . \\%
                \global\no@bbfalse
        \fi
   \fi\fi}%
%
%
\def\epsf@atendlit{(atend)} 
\def\epsf@atend #1 #2 #3\\{%
   \def\epsf@tmp{#1}\ifx\epsf@tmp\empty
      \epsf@atend #2 #3 .\\\else
   \ifx\epsf@tmp\epsf@atendlit\@atendtrue\fi\fi}


\chardef\letter = 11
\chardef\other = 12

\newif \ifdebug 
\newif\ifc@mpute 
\newif\if@atend
\c@mputetrue 

\let\then = \relax
\def\r@dian{pt }
\let\r@dians = \r@dian
\let\dimensionless@nit = \r@dian
\let\dimensionless@nits = \dimensionless@nit
\def\internal@nit{sp }
\let\internal@nits = \internal@nit
\newif\ifstillc@nverging
\def \Mess@ge #1{\ifdebug \then \message {#1} \fi}

{ 
        \catcode `\@ = \letter
        \gdef \nodimen {\expandafter \n@dimen \the \dimen}
        \gdef \term #1 #2 #3%
               {\edef \t@ {\the #1}
                \edef \t@@ {\expandafter \n@dimen \the #2\r@dian}%
                \t@rm {\t@} {\t@@} {#3}%
               }
        \gdef \t@rm #1 #2 #3%
               {{%
                \count 0 = 0
                \dimen 0 = 1 \dimensionless@nit
                \dimen 2 = #2\relax
                \Mess@ge {Calculating term #1 of \nodimen 2}%
                \loop
                \ifnum  \count 0 < #1
                \then   \advance \count 0 by 1
                        \Mess@ge {Iteration \the \count 0 \space}%
                        \Multiply \dimen 0 by {\dimen 2}%
                        \Mess@ge {After multiplication, term = \nodimen 0}%
                        \Divide \dimen 0 by {\count 0}%
                        \Mess@ge {After division, term = \nodimen 0}%
                \repeat
                \Mess@ge {Final value for term #1 of 
                                \nodimen 2 \space is \nodimen 0}%
                \xdef \Term {#3 = \nodimen 0 \r@dians}%
                \aftergroup \Term
               }}
        \catcode `\p = \other
        \catcode `\t = \other
        \gdef \n@dimen #1pt{#1} 
}

\def \Divide #1by #2{\divide #1 by #2} 

\def \Multiply #1by #2
       {{
        \count 0 = #1\relax
        \count 2 = #2\relax
        \count 4 = 65536
        \Mess@ge {Before scaling, count 0 = \the \count 0 \space and
                        count 2 = \the \count 2}%
        \ifnum  \count 0 > 32767 
        \then   \divide \count 0 by 4
                \divide \count 4 by 4
        \else   \ifnum  \count 0 < -32767
                \then   \divide \count 0 by 4
                        \divide \count 4 by 4
                \else
                \fi
        \fi
        \ifnum  \count 2 > 32767 
        \then   \divide \count 2 by 4
                \divide \count 4 by 4
        \else   \ifnum  \count 2 < -32767
                \then   \divide \count 2 by 4
                        \divide \count 4 by 4
                \else
                \fi
        \fi
        \multiply \count 0 by \count 2
        \divide \count 0 by \count 4
        \xdef \product {#1 = \the \count 0 \internal@nits}%
        \aftergroup \product
       }}

\def\r@duce{\ifdim\dimen0 > 90\r@dian \then   
                \multiply\dimen0 by -1
                \advance\dimen0 by 180\r@dian
                \r@duce
            \else \ifdim\dimen0 < -90\r@dian \then  
                \advance\dimen0 by 360\r@dian
                \r@duce
                \fi
            \fi}

\def\Sine#1%
       {{%
        \dimen 0 = #1 \r@dian
        \r@duce
        \ifdim\dimen0 = -90\r@dian \then
           \dimen4 = -1\r@dian
           \c@mputefalse
        \fi
        \ifdim\dimen0 = 90\r@dian \then
           \dimen4 = 1\r@dian
           \c@mputefalse
        \fi
        \ifdim\dimen0 = 0\r@dian \then
           \dimen4 = 0\r@dian
           \c@mputefalse
        \fi
        \ifc@mpute \then
                \divide\dimen0 by 180
                \dimen0=3.141592654\dimen0
                \dimen 2 = 3.1415926535897963\r@dian 
                \divide\dimen 2 by 2 
                \Mess@ge {Sin: calculating Sin of \nodimen 0}%
                \count 0 = 1 
                \dimen 2 = 1 \r@dian 
                \dimen 4 = 0 \r@dian 
                \loop
                        \ifnum  \dimen 2 = 0 
                        \then   \stillc@nvergingfalse 
                        \else   \stillc@nvergingtrue
                        \fi
                        \ifstillc@nverging 
                        \then   \term {\count 0} {\dimen 0} {\dimen 2}%
                                \advance \count 0 by 2
                                \count 2 = \count 0
                                \divide \count 2 by 2
                                \ifodd  \count 2 
                                \then   \advance \dimen 4 by \dimen 2
                                \else   \advance \dimen 4 by -\dimen 2
                                \fi
                \repeat
        \fi             
                        \xdef \sine {\nodimen 4}%
       }}

\def\Cosine#1{\ifx\sine\UnDefined\edef\Savesine{\relax}\else
                             \edef\Savesine{\sine}\fi
        {\dimen0=#1\r@dian\multiply\dimen0 by -1
         \advance\dimen0 by 90\r@dian
         \Sine{\nodimen 0}
         \xdef\cosine{\sine}
         \xdef\sine{\Savesine}}}              
%
\def\psdraft{\def\@psdraft{0}}
\def\psfull{\def\@psdraft{100}}
\psfull
\newif\if@draftbox
\def\psnodraftbox{\@draftboxfalse}
\@draftboxtrue
\newif\if@noisy
\def\pssilent{\@noisyfalse}
\def\psnoisy{\@noisytrue}
\@noisyfalse

\newif\if@bbllx
\newif\if@bblly
\newif\if@bburx
\newif\if@bbury
\newif\if@height
\newif\if@width
\newif\if@rheight
\newif\if@rwidth
\newif\if@angle
\newif\if@clip
\newif\if@verbose
\def\@p@@sclip#1{\@cliptrue}
\newif\if@prologfile
\@prologfiletrue
%
\newif\ifuse@psfig
\def\@p@@sfile#1{%
\def\@p@sfile{NO FILE: #1}%
\def\@p@sfilefinal{NO FILE: #1}%
        \openin1=#1
        \ifeof1\closein1
                \openin1=#1.bb
                        \ifeof1\closein1
                                \if@bbllx\if@bblly\if@bburx\if@bbury
                                        \def\@p@sfile{#1}%
                                        \def\@p@sfilefinal{#1}%
                                        \fi\fi\fi
                                \else
                                        \@latexerr{ERROR! PostScript file #1 not found}\@whattodo
                                        \@p@@sbbllx{100bp}
                                        \@p@@sbblly{100bp}
                                        \@p@@sbburx{200bp}
                                        \@p@@sbbury{200bp}
                                        \def\@p@scost{200}
                                \fi
                        \else
                                \closein1%
                                \edef\@p@sfile{#1.bb}%
                                \edef\@p@sfilefinal{"`zcat `texfind #1.Z`"}%
                        \fi
        \else\closein1
                    \edef\@p@sfile{#1}%
                    \edef\@p@sfilefinal{#1}%
        \fi%
}
\let\@p@@sfigure\@p@@sfile
\def\@p@@sbbllx#1{
                \use@psfigtrue
                \@bbllxtrue
                \dimen100=#1
                \edef\@p@sbbllx{\number\dimen100}
}
\def\@p@@sbblly#1{
                \use@psfigtrue
                \@bbllytrue
                \dimen100=#1
                \edef\@p@sbblly{\number\dimen100}
}
\def\@p@@sbburx#1{
                \use@psfigtrue
                \@bburxtrue
                \dimen100=#1
                \edef\@p@sbburx{\number\dimen100}
}
\def\@p@@sbbury#1{
                \use@psfigtrue
                \@bburytrue
                \dimen100=#1
                \edef\@p@sbbury{\number\dimen100}
}
\def\@p@@sheight#1{
                \@heighttrue
                \epsfysize=#1
                \dimen100=#1
                \edef\@p@sheight{\number\dimen100}
}
\def\@p@@swidth#1{
                \@widthtrue
                \epsfxsize=#1
                \dimen100=#1
                \edef\@p@swidth{\number\dimen100}
}
\def\@p@@srheight#1{
                \@rheighttrue
                \dimen100=#1
                \edef\@p@srheight{\number\dimen100}
}
\def\@p@@srwidth#1{
                \@rwidthtrue
                \dimen100=#1
                \edef\@p@srwidth{\number\dimen100}
}
\def\@p@@sangle#1{
                \use@psfigtrue
                \@angletrue
                \edef\@p@sangle{#1} 
}
\def\@p@@ssilent#1{ 
                \@verbosefalse
}
\def\@cs@name#1{\csname #1\endcsname}
\def\@setparms#1=#2,{\@cs@name{@p@@s#1}{#2}}
%
%
\def\ps@init@parms{
                \@bbllxfalse \@bbllyfalse
                \@bburxfalse \@bburyfalse
                \@heightfalse \@widthfalse
                \@rheightfalse \@rwidthfalse
                \def\@p@sbbllx{}\def\@p@sbblly{}
                \def\@p@sbburx{}\def\@p@sbbury{}
                \def\@p@sheight{}\def\@p@swidth{}
                \def\@p@srheight{}\def\@p@srwidth{}
                \def\@p@sangle{0}
                \def\@p@sfile{}
                \def\@p@scost{10}
                \use@psfigfalse
                \def\@sc{}
                \if@noisy
                        \@verbosetrue
                \else
                        \@verbosefalse
                \fi
                \@clipfalse
}
%
%
\def\parse@ps@parms#1{
                \@psdo\@psfiga:=#1\do
                   {\expandafter\@setparms\@psfiga,}}
%
%
\newif\ifno@bb
\def\bb@missing{
        \epsfgetbb{\@p@sfile}
        \ifepsfbbfound\no@bbfalse\else\no@bbtrue\bb@cull\epsfllx\epsflly\epsfurx\epsfury\fi
}       
\def\bb@cull#1#2#3#4{
        \dimen100=#1 bp\edef\@p@sbbllx{\number\dimen100}
        \dimen100=#2 bp\edef\@p@sbblly{\number\dimen100}
        \dimen100=#3 bp\edef\@p@sbburx{\number\dimen100}
        \dimen100=#4 bp\edef\@p@sbbury{\number\dimen100}
        \no@bbfalse
}

\newdimen\p@intvaluex
\newdimen\p@intvaluey
\newdimen\@ffsetvalue
\newdimen\x@ffsetvalue
\newdimen\y@ffsetvalue


\def\compute@offset#1#2{{\dimen0=#1 sp\dimen1=#2 sp
                        \advance\dimen1 by -\dimen0
                        \dimen1=\sine\dimen1
                        \dimen0=\cosine\dimen1
                        \ifdim\dimen0<0sp \dimen1=0sp \fi
                        \global\@ffsetvalue=\dimen1}}

\def\rotate@#1#2{{\dimen0=#1 sp\dimen1=#2 sp
                  \global\p@intvaluex=\cosine\dimen0
                  \dimen3=\sine\dimen1
                  \global\advance\p@intvaluex by -\dimen3
                  \global\p@intvaluey=\sine\dimen0
                  \dimen3=\cosine\dimen1
                  \global\advance\p@intvaluey by \dimen3
                  }}
%
\def\compute@bb{
                \no@bbfalse
                \if@bbllx \else \no@bbtrue \fi
                \if@bblly \else \no@bbtrue \fi
                \if@bburx \else \no@bbtrue \fi
                \if@bbury \else \no@bbtrue \fi
                \ifno@bb \bb@missing \fi
                \ifno@bb
                        \@latexerr{ERROR! cannot locate BB!}\@whattodobb
                        \@p@@sbbllx{100bp}
                        \@p@@sbblly{100bp}
                        \@p@@sbburx{200bp}
                        \@p@@sbbury{200bp}
                        \def\@p@scost{200}
                \fi
                \if@angle 
                        \Sine{\@p@sangle}\Cosine{\@p@sangle}
                        \compute@offset{\@p@sbblly}{\@p@sbbury}
                        \x@ffsetvalue=\@ffsetvalue
                        \compute@offset{\@p@sbburx}{\@p@sbbllx}
                        \y@ffsetvalue=\@ffsetvalue

                        \rotate@{\@p@sbbllx}{\@p@sbblly}
                        \advance\p@intvaluex by -\x@ffsetvalue
                        \advance\p@intvaluey by -\y@ffsetvalue
                        \edef\@p@sbbllx{\number\p@intvaluex}
                        \edef\@p@sbblly{\number\p@intvaluey}

                        \rotate@{\@p@sbburx}{\@p@sbbury}
                        \advance\p@intvaluex by \x@ffsetvalue
                        \advance\p@intvaluey by \y@ffsetvalue
                        \edef\@p@sbburx{\number\p@intvaluex}
                        \edef\@p@sbbury{\number\p@intvaluey}
                        {
                         \count0=\@p@sbbllx \count1=\@p@sbblly
                         \count2=\@p@sbburx \count3=\@p@sbbury
                         \dimen0=\@p@sbbllx sp\dimen1=\@p@sbblly sp
                         \dimen2=\@p@sbburx sp\dimen3=\@p@sbbury sp
                         \dimen203=\dimen2 \advance\dimen203 by -\dimen0
                         \dimen204=\dimen3 \advance\dimen204 by -\dimen1
                         \ifdim\dimen203<0sp 
                              \count203=\count2 \count2=\count0 
                              \count0=\count203 
                              \global\edef\@p@sbbllx{\number\count0}
                              \global\edef\@p@sbburx{\number\count2}
                         \fi
                         \ifdim\dimen204<0sp 
                               \count204=\count3
                               \count3=\count1
                               \count1=\count204
                               \global\edef\@p@sbblly{\number\count1}
                               \global\edef\@p@sbbury{\number\count3}
                         \fi
                        }
                \fi
                \count203=\@p@sbburx
                \count204=\@p@sbbury
                \advance\count203 by -\@p@sbbllx
                \advance\count204 by -\@p@sbblly
                \edef\@bbw{\number\count203}
                \edef\@bbh{\number\count204}
}
%
%
\def\in@hundreds#1#2#3{\count240=#2 \count241=#3
                     \count100=\count240        
                     \divide\count100 by \count241
                     \count101=\count100
                     \multiply\count101 by \count241
                     \advance\count240 by -\count101
                     \multiply\count240 by 10
                     \count101=\count240        
                     \divide\count101 by \count241
                     \count102=\count101
                     \multiply\count102 by \count241
                     \advance\count240 by -\count102
                     \multiply\count240 by 10
                     \count102=\count240        
                     \divide\count102 by \count241
                     \count200=#1\count205=0
                     \count201=\count200
                        \multiply\count201 by \count100
                        \advance\count205 by \count201
                     \count201=\count200
                        \divide\count201 by 10
                        \multiply\count201 by \count101
                        \advance\count205 by \count201
                     \count201=\count200
                        \divide\count201 by 100
                        \multiply\count201 by \count102
                        \advance\count205 by \count201
                     \edef\@result{\number\count205}
}
\def\compute@wfromh{
                \in@hundreds{\@p@sheight}{\@bbw}{\@bbh}
                \edef\@p@swidth{\@result}
}
\def\compute@hfromw{
                \in@hundreds{\@p@swidth}{\@bbh}{\@bbw}
                \edef\@p@sheight{\@result}
}
\def\compute@handw{
                \if@height 
                        \if@width
                        \else
                                \compute@wfromh
                        \fi
                \else 
                        \if@width
                                \compute@hfromw
                        \else
                                \edef\@p@sheight{\@bbh}
                                \edef\@p@swidth{\@bbw}
                        \fi
                \fi
}
\def\compute@resv{
                \if@rheight \else \edef\@p@srheight{\@p@sheight} \fi
                \if@rwidth \else \edef\@p@srwidth{\@p@swidth} \fi
}
%
\def\compute@sizes{
        \compute@bb
        \compute@handw
        \compute@resv
}
%
%
\def\psfig#1{\vbox {
        %
        \ps@init@parms
        \parse@ps@parms{#1}
        \ifnum\@p@scost<\@psdraft
                \typeout{[\@p@sfilefinal]}
                \if@verbose
                        \typeout{epsfig: using PSFIG macros}
                \fi
                \psfig@method
        \else
                \epsfig@draft
        \fi
}}

\def\epsfig#1{\vbox {
        %
        \ps@init@parms
        \parse@ps@parms{#1}
        \ifnum\@p@scost<\@psdraft
                \typeout{[\@p@sfilefinal]}
                \if@clip\use@psfigtrue\fi
                \if@angle\use@psfigtrue\fi
                \ifuse@psfig
                        \if@verbose
                                \typeout{epsfig: using PSFIG macros}
                        \fi
                        \psfig@method
                \else
                        \if@verbose
                                \typeout{epsfig: using EPSF macros}
                        \fi
                        \epsf@method
                \fi
        \else
                \epsfig@draft
        \fi
}}

\def\epsf@method{%
        \epsfgetbb{\@p@sfile}%
        \epsfsetgraph{\@p@sfilefinal}
}
\def\psfig@method{%
        \compute@sizes
        \special{ps::[begin]  \@p@swidth \space \@p@sheight \space%
        \@p@sbbllx \space \@p@sbblly \space%
        \@p@sbburx \space \@p@sbbury \space%
        startTexFig \space }%
        \if@angle
                \special {ps:: \@p@sangle \space rotate \space} 
        \fi
        \if@clip
                \if@verbose
                        \typeout{(clipped to BB) }
                \fi
                \special{ps:: doclip \space }%
        \fi
        \special{ps: plotfile \@p@sfilefinal \space }%
        \special{ps::[end] endTexFig \space }%
        \vbox to \@p@srheight true sp{\hbox to \@p@srwidth true sp{\hss}\vss}
}
%
\def\epsfig@draft{
\compute@sizes
\if@draftbox
        \hbox{\fbox{\vbox to \@p@srheight true sp{
        \vss\hbox to \@p@srwidth true sp{ \hss \@p@sfilefinal \hss }\vss
        }}}
\else
        \vbox to \@p@srheight true sp{%
        \vss\hbox to \@p@srwidth true sp{\hss}\vss}
\fi     
}
\catcode`\@=\atcode 

\def\tilde{$\sim$}

\newcount\dcoef
\newcount\mcoef

\def\dessin(#1)(#2)(#3){{\epsfig{file=#1,width=#2truecm,height=#3truecm}}}

\mcoef=1
\dcoef=2